# LIMIT THEOREMS FOR MIXED MAX–SUM PROCESSES WITH RENEWAL STOPPING

By Dmitrii S. Silvestrov[1] and Jozef L. Teugels

*Mälardalen University and Katholieke Universiteit Leuven*

This article is devoted to the investigation of limit theorems for mixed max–sum processes with renewal type stopping indexes. Limit theorems of weak convergence type are obtained as well as functional limit theorems.

**1. Introduction.** The main object of this article is the derivation of a number of limit theorems for mixed max–sum processes with renewal stopping. Such processes are constructed from the three-component sequences of i.i.d. random vectors taking values in $R_1 \times R_1 \times [0, \infty)$ in the following way. The first component of the sequence is used to construct an extremal max process of i.i.d. random variables. The second is used as a traditional real-valued sum process of i.i.d. random variables. Finally, the third component is introduced by a nonnegative sum process of i.i.d. random variables. It induces the stopping renewal process that is a process of the first exceeding times over a specific level $t > 0$. The first two components are then stopped using this renewal process. The overall process so obtained will be called a *max–sum process with renewal stopping*. Note that at this point we do not restrict possible dependencies between the three components.

Max–sum processes with renewal stopping of the above type naturally appear in various applications. To help visualize such processes, we give a few concrete examples.

Example 1. Consider an ordinary renewal process $\{X_1, X_2, \ldots\}$ generated by nonnegative independent random variables with common distribution $F$ of $X$. Define the *renewal counting process* $N(t) = \max(n : \sum_{i \le n} X_i \le$

---

Received September 2001; revised June 2003.

[1]Supported in part by Grant F/98/093 from the Catholic University of Leuven.

---

$t$). The two-dimensional process $(\max_{i \leq N(t)} X_i, \sum_{i \leq N(t)} X_i)$ serves as a special case of a max–sum process with renewal stopping where max, sum and renewal components are constructed from the same $X$ sequence of i.i.d. random variables.

EXAMPLE 2. Consider a similar setup based on a sequence $\{(X_1, Y_1), (X_2, Y_2), \ldots\}$ of i.i.d. random vectors with nonnegative components. Interpret the $X$ values as the times in between claims in an insurance portfolio and $Y$ as the claim values. The renewal counting process $N(t)$, constructed as in the first example, counts the number of claims in the time interval $[0, t]$. The two-dimensional process $(\max_{i \leq N(t)} Y_i, \sum_{i \leq N(t)} Y_i)$ serves as a special case of a max–sum process with renewal stopping, where the first two components are constructed from the same $Y$ sequence of i.i.d. random variables, while the renewal stopping process is constructed from the $X$ sequence. The first component represents the maximal claim experienced within the portfolio over the time interval $[0, t]$, while the second component models the totality of claims over the same period. The possible dependence between arrival times and claim sizes can be kept as part of the model.

EXAMPLE 3. Consider once more a similar setup, but now based on a sequence of triplets $\{(X_1, Y_1, Z_1), (X_2, Y_2, Z_2), \ldots\}$ of i.i.d. random vectors with nonnegative components. Interpret the $X$ values as the interarrival times between earthquakes at a specific location, the $Y$ values as the sizes of the quakes and the $Z$ values as the damages caused by these quakes. Alternatively, the last component could also model the number of aftershocks from the corresponding earthquake. The renewal counting process $N(t)$, constructed as in the first example, counts the number of quakes in the time interval $[0, t]$. Again, the bivariate process $(\max_{i \leq N(t)} Y_i, \sum_{i \leq N(t)} Z_i)$ illustrates the max-sum process as constructed from a three-dimensional $(X, Y, Z)$ sequence. It is rather natural to look at the maximal size of the quakes, represented by the first max component. The second sum component could refer to the total damage caused by it. In the alternative interpretation, the last component models the number of aftershocks from this major earthquake.

In all three examples, the question of interest is the influence of the first max component on the second sum component. As a particular case, one might be interested in the asymptotic behavior when time $t$ tends to infinity. The answer to this question relies on another problem that needs to be solved first. What is the asymptotic behavior of the joint distribution of the corresponding mixed max–sum processes with renewal stopping? This is precisely the subject of the present article.

If we look only at the first and third component, then we are in the realm of limit theory for extremal processes with random sample size index. This area has been thoroughly studied by Berman (1962), Barndorff-Nielsen (1964), Mogyoródi (1967), Thomas (1972), Sen (1972), Galambos (1973, 1975, 1978, 1992, 1994), Gnedenko and Gnedenko (1982), Beirlant and Teugels (1992) and, more recently, Silvestrov and Teugels (1998a, b). We need to point out that the last article covers limit theorems for extremal max processes with renewal stopping for the case of asymptotic independency between the extremal process and the renewal stopping process. Some related results concerning exceedances of ergodic regenerative processes with discrete time can also be found in works by Serfozo (1980), Rootzén (1988) and Leadbetter and Rootzén (1988).

Results concerning the joint asymptotic behavior of maxima and sums of i.i.d. random variables, in particular concerning the case when their quotient tends to 1, were obtained by Arov and Bobrov (1960), O'Brien (1980), Maller and Resnick (1984) and Pruitt (1987). Related results can also be found in Darling (1952), Smirnov (1952) and Aebi, Embrechts and Mikosch (1992). In the general case, the joint asymptotic distribution of maxima and sums of i.i.d. random variables was studied for the scale–location model by Breiman (1965), Chow and Teugels (1979), Resnick (1986) and Haas (1992). Related results also can be found in Lamperti (1964), Anderson and Turkman (1991), Kesten and Maller (1994), Hsing (1995), Ho and Hsing (1996) and the book edited by Hahn, Mason and Weiner (1991).

Looking only at the second and third components, we end up in the area of limit theorems for sums of random variables with random index. To avoid overloading the bibliography, we refer only to contributions devoted to theorems for sums with renewal type stopping indexes. A good bibliography for the period up to the early seventies is Serfozo (1975). Another good summary can be found in the book by Gut (1988). More specific references are Feller (1949, 1966), Smith (1955), Dynkin (1955), Takacs (1959), Lamperti (1961, 1962), Borovkov (1967), Iglehart (1969), Silvestrov (1972a, 1974, 1983, 1991), Kaplan and Silvestrov (1979), Gut and Janson (1983), Niculescu (1984), Murphree and Smith (1986), Shedler (1987, 1993) and Roginsky (1989, 1994). Another source of results is formed by diffusion approximations for risk processes as in Iglehart (1969), Siegmund (1975), Harrison (1977), Gerber (1979), Grandell (1977, 1991), Beard, Pentikäinen and Pesonen (1984), Aebi, Embrechts and Mikosch (1994), Asmussen (1984, 2000, 2003), Schmidli (1992, 1997) and the books by Embrechts, Klüppelberg and Mikosch (1997), Rolski, Schmidli, Schmidt and Teugels (1999) and Bening and Korolev (2002).

We also refer to the articles by Silvestrov and Teugels (2001), which is an extended report version of the current article, and Silvestrov (2000, 2002), where one can find an extended bibliography of publications from the realm of this article.

The goal of the current work is to derive weak and functional limit theorems for a combination of the above processes that we have coined max–sum processes with renewal stopping. Our model includes as particular cases all three types of models mentioned above, that is, mixed max–sum processes, max processes with random indices of renewal type and sum processes with renewal stopping. Finding inspiration in the classical model with i.i.d. random variables, the model under consideration deals nevertheless with the joint behavior of max processes, sum processes and renewal stopping processes where dependencies can be introduced via the components of the initial i.i.d. random vectors. To increase their generality and applicability, all results are presented in a random process setting and for a general triangular array model. We show that the corresponding limit theorems can be obtained under conditions analogous to the well-known conditions of the central criterion of convergence for sums of i.i.d. random variables and similar conditions for maxima of i.i.d. random variables. No additional technical assumptions are involved. In this sense, limit theorems presented herein give some kind of a "final solution" for limit theorems for max–sum processes with renewal stopping based on i.i.d. random variables.

More specifically, let $\{(\xi_{\varepsilon,n}, \gamma_{\varepsilon,n}, \kappa_{\varepsilon,n}), \ n = 1, 2, \ldots\}$ be for every $\varepsilon > 0$ a sequence of i.i.d. random variables taking values in $R_1 \times R_1 \times [0, \infty)$. Furthermore, we need nonrandom functions $n_\varepsilon > 0$ for which $n_\varepsilon \to \infty$ as $\varepsilon \to 0$.

We first introduce the components using nonrandom sample sizes, namely the *extremal max process* $\xi_\varepsilon(t) = \max_{k \leq 1 \vee tn_\varepsilon} \xi_{\varepsilon,k}, \ t \geq 0$, the *sum process* $\gamma_\varepsilon(t) = \sum_{k \leq tn_\varepsilon} \gamma_{\varepsilon,k}, \ t \geq 0$, and the *positive sum process* $\kappa_\varepsilon(t) = \sum_{k \leq tn_\varepsilon} \kappa_{\varepsilon,k}, \ t \geq 0$. The last process induces a *renewal stopping process* $\tau_\varepsilon(t) = \sup(s : \kappa_\varepsilon(s) \leq t), \ t > 0$, that will behave like a random time clock on the three separate components.

In the current article our attention goes to a thorough study of the *max–sum process with renewal stopping* $(\xi_\varepsilon(\tau_\varepsilon(t)), \gamma_\varepsilon(\tau_\varepsilon(t)), \kappa_\varepsilon(\tau_\varepsilon(t))), \ t > 0$, because we will give general conditions for weak convergence of these processes and for their functional counterparts.

Before doing that we need to investigate in detail limit theorems for the three-dimensional mixed max–sum processes $\{(\xi_\varepsilon(t), \gamma_\varepsilon(t), \kappa_\varepsilon(t)), \ t > 0\}$. Among others, we give conditions of weak convergence of such processes as well as conditions of their convergence in Skorokhod $J$-topology. We therefore first recall what is known about the separate components of this three-variate process. Because we definitely need marginal weak convergence of the max processes $\{\xi_\varepsilon(t), t > 0\}$ and of the two-component sum processes $\{(\kappa_\varepsilon(t), \gamma_\varepsilon(t)), t \geq 0\}$, we automatically need a set of necessary conditions. They help in the formulation of the conditions for joint weak convergence, in particular for the corresponding limiting three-component mixed max–sum process.

The weak convergence and $J$-convergence of mixed max–sum processes is treated in Section 2. In Section 3 we deal with the weak and the $J$-convergence of the max–sum processes with renewal stopping. In the final Section 4 we treat some examples.

**2. Mixed max–sum processes.** In this section we deal with general conditions for the weak and functional convergence of the mixed max–sum processes. We start with the extremal component. We then turn to the sum processes to finish this section with the mixed max–sum processes.

2.1. *Weak convergence of max processes.* We start with the extremal component. As usual, let us denote $C_f$ the set of continuity of a function $f$. The following condition is standard in papers dealing with limit theorems for extremes:

CONDITION A. As $\varepsilon > 0$, $n_\varepsilon P\{\xi_{\varepsilon,1} > u\} \to \pi_1(u)$ for all $u \in R_1$ which belongs to the set $C_{\pi_1}$ for the limiting function $\pi_1(u)$.

The ingredients on the right satisfy a number of conditions.

- The function $\pi_1(u)$ acts from $(-\infty, \infty)$ into $[0, \infty]$ and is nonincreasing and continuous from the right: if $\pi_1(u) = \infty$, continuity is interpreted as $\pi_1(t) \uparrow \infty$ as $t \downarrow u$; furthermore, $\pi_1(-\infty) = \infty$ and $\pi_1(\infty) = 0$.
- As such, these conditions imply that the function $\exp(-\pi_1(u))$ is a distribution function. If we define $u_\pi = \sup(u : \pi_1(u) = \infty) \geq -\infty$, then $\exp(-\pi_1(u))$ takes positive values for $u > u_\pi$ and $\exp(-\pi_1(u)) = 0$ for $u < u_\pi$, while $\exp(-\pi_1(u_\pi)) = \exp(-\pi_1(u_\pi+))$ can take any value in the interval $[0, 1]$.

One of the important aspects of classical extreme value theory is the *scale–location model*. Here the random variables $\xi_{\varepsilon,n}$ are represented in the form $\xi_{\varepsilon,n} = (\xi_n - a_\varepsilon)/b_\varepsilon$, where $\xi_n, n = 1, 2, \ldots$, are i.i.d. random variables, and $a_\varepsilon$ and $b_\varepsilon$ are some nonrandom centralization and normalization constants. In this case, the distribution $\exp(-\pi_1(u))$ belongs to one of three families of classical extremal distributions. See, for instance, Galambos (1978), Leadbetter, Lindgren and Rootzén (1983), Resnick (1987) and Berman (1992).

This one-dimensional result can be extended. Denote by $D_0$ the space of step functions on $(0, \infty)$ continuous from the right and with a finite number of only positive jumps in every finite subinterval of $(0, \infty)$. It is known [see, e.g., Serfozo (1982), Leadbetter, Lindgren and Rootzén (1983), Resnick (1987) and Berman (1992)] that Condition A is necessary and sufficient for the weak convergence

$$(2.1) \qquad \{\xi_\varepsilon(t), t > 0\} \quad \Longrightarrow \quad \{\xi_0(t), t > 0\} \qquad \text{as } \varepsilon \to 0.$$

The limiting process $\{\xi_0(t), t > 0\}$ in (2.1) is called an *extremal process*. It has the following finite-dimensional distributions for $0 = t_0 < t_1 < \cdots < t_n$, $-\infty < u_1 \leq \cdots \leq u_n < \infty$, $n \geq 1$:

$$(2.2) \quad P\{\xi_0(t_1) \leq u_1, \ldots, \xi_0(t_n) \leq u_n\} = \prod_{k=1}^{n} \exp(-\pi_1(u_k)(t_k - t_{k-1})).$$

Referring to the notation $u_\pi$ above, let $v_\pi = \inf(u : \pi_1(u) = 0) \leq \infty$. Then the distribution function $\exp(-\pi_1(u))$ is concentrated on the interval $[u_\pi, v_\pi]$, $\xi_0(t) \to u_\pi$ a.s. as $t \to 0$, while $\xi_0(t) \to v_\pi$ a.s. as $t \to \infty$.

Also, the process $\{\xi_0(t), t > 0\}$ is a stochastically continuous homogeneous Markov jump process whose trajectories belong to the space $D_0$ with probability 1. It has transition probabilities

$$(2.3) \quad P\{\xi_0(s + t) \leq u | \xi_0(s) = v\} = \chi(v \leq u) \exp(-t\pi_1(u)),$$

where $\chi(A)$ is used for the indicator of event $A$.

2.2. *Weak convergence of sum processes.* Let us consider the bivariate process $\{(\gamma_\varepsilon(t), \kappa_\varepsilon(t)), t > 0\}$. It is a process of step sums of i.i.d. random vectors. Conditions of weak convergence of these processes can be formulated with the use of the vector form of the classical criterion for weak convergence given in Skorokhod (1964). These conditions involve the tail probabilities, the truncated means and the truncated variances of the random variables.

We use the abbreviation $R = (R_1 \times [0, \infty)) \setminus \{(0, 0)\}$ and we write $B_R$ for the Borel $\sigma$-algebra of subsets of $R$. Let $\Phi$ be the class of continuous bounded functions defined on $R$ and vanishing in some neighborhood of the point $(0, 0)$. The condition for convergence takes the following form:

CONDITION B.  (a) As $\varepsilon \to 0$, $n_\varepsilon E\phi(\gamma_{\varepsilon,1}, \kappa_{\varepsilon,1}) = n_\varepsilon \int_R \phi(v, w) P\{(\gamma_{\varepsilon,1}, \kappa_{\varepsilon,1}) \in dv \times dw\} \to \int_R \phi(v, w)\Pi_{2,3}(dv \times dw)$ for all $\phi \in \Phi$.

(b) As $\varepsilon \to 0$, $n_\varepsilon E\gamma_{\varepsilon,1}\chi(|\gamma_{\varepsilon,1}| \leq v) \to a(v)$ for some $v > 0$ for which the points $\pm v$ are points of continuity of the function $\Pi_{2,3}(\{v\} \times [0, \infty))$.

(c) As $\varepsilon \to 0$, $n_\varepsilon E\kappa_{\varepsilon,1}\chi(\kappa_{\varepsilon,1} \leq w) \to c(v)$ for some $w > 0$ which is a point of continuity of the function $\Pi_{2,3}(R_1 \times \{w\})$.

(d) As $\varepsilon \to 0$ and then $v \to 0$, $n_\varepsilon(E\gamma_{\varepsilon,1}^2\chi(|\gamma_{\varepsilon,1}| \leq v) - (E\gamma_{\varepsilon,1}\chi(|\gamma_{\varepsilon,1}| \leq v))^2) \to b^2$. This expression refers to two repeated limits of the form $\lim_{0 < v \to 0} \limsup_{\varepsilon \to 0}$ and $\lim_{0 < v \to 0} \liminf_{\varepsilon \to 0}$.

We list a number of properties of the limits on the right.

- $\Pi_{2,3}(A)$ is a measure on the $\sigma$-algebra $B_R$.
- The projection $\Pi_2(B) = \Pi_{2,3}(B \times [0, \infty))$ is a measure on the Borel $\sigma$-algebra of subsets of $R_1 \setminus \{0\}$ such that $\fint_{R_1} s^2/(s^2 + 1)\Pi_2(ds) < \infty$, where $\fint$ is an integral over the corresponding interval with the point 0 excluded from the interval of integration.

- The projection $\Pi_3(C) = \Pi_{2,3}(R_1 \times C)$ is a measure on the Borel $\sigma$-algebra of subsets of $(0, \infty)$ such that $\int_{(0,\infty)} s/(s+1)\Pi_3(ds) < \infty$.

- $a(v), v > 0$ is a real-valued measurable function and Condition B(b) can, under Condition B(a), only hold simultaneously for all points $v > 0$ for which $\pm v$ are points of continuity of the function $\Pi_2(\{v\})$ and for any such point the constant $a = a(v) - \oint_{|s|<v} s^3/(1+s^2)\Pi_2(ds) + \int_{|s|>v} s/(1+s^2)\Pi_2(ds)$ does not depend on the choice of $v$.

- Function $c(w), w > 0$, is a nonnegative nondecreasing function and Condition B(c) can, under Condition B(a), only hold simultaneously for all points $w > 0$ for which $w$ is a point of continuity of the function $\Pi_3(\{w\})$ and for any such point the constants $c = c(w) - \int_{(0,w)} s\Pi_3(ds) \geq 0$ and $d = c + \int_{(0,\infty)} s/(1+s^2)\Pi_3(ds)$ do not depend on the choice of $w$.

- Finally, $b^2$ is a nonnegative constant.

The nonnegativity of the random variables $\kappa_{\varepsilon,1}$ and Conditions B(a)–(c) imply that the repeated limits for variances of these random variables and covariances of random variables $\kappa_{\varepsilon,1}$ and $\gamma_{\varepsilon,1}$, analogous to those in Condition B(d), are equal to zero. For this reason, the corresponding conditions are not included in Condition B.

According to a central criterion of convergence [see, e.g., Loève (1955)], Condition B is necessary and sufficient for the bivariate weak convergence

$$(2.4) \quad \{(\gamma_\varepsilon(t), \kappa_\varepsilon(t)), t \geq 0\} \quad \Longrightarrow \quad \{(\gamma_0(t), \kappa_0(t)), t \geq 0\} \qquad \text{as } \varepsilon \to 0.$$

Denote by $D$ the space of functions on $(0, \infty)$ without discontinuities of the second kind and continuous from the right. Let also $D_+$ be the space of nondecreasing functions from $D$ and let $D_{++}$ be the space of nonnegative functions from $D_+$.

The limiting process $\{(\gamma_0(t), \kappa_0(t)), t > 0\}$ is a homogeneous stochastically continuous process with independent increments whose trajectories belong to the space $D \times D_{++}$ with probability 1 and with the characteristic function for $t > 0$ given by

$$
(2.5) \quad
\begin{aligned}
E\exp&\{i(y\gamma_0(t) + z\kappa_0(t))\}\\
&= \varphi_{2,3}(t, y, z)\\
&= \exp\Big\{ t\Big( iay - \frac{1}{2}b^2y^2 + i\,dz \\
&\qquad + \int_R \Big( e^{i(yv+zw)} - 1 - \frac{i(yv+zw)}{1+v^2+w^2}\Big)\Pi_{2,3}(dv \times dw)\Big)\Big\}.
\end{aligned}
$$

2.3. *Weak convergence of mixed max–sum processes.* We finally turn to the study of the joint behavior of the three components together. The following condition should be added to Conditions A and B to provide joint weak convergence of max and sum processes:

CONDITION C. As $\varepsilon \to 0$, $n_\varepsilon E\chi(\xi_{\varepsilon,1} > u)\phi(\gamma_{\varepsilon,1}, \kappa_{\varepsilon,1}) = n_\varepsilon \times \int_R \phi(v,w)P\{\xi_{\varepsilon,1} > u, (\gamma_{\varepsilon,1}, \kappa_{\varepsilon,1}) \in dv \times dw\} \to \int_R \phi(v,w)\Pi_{2,3}^{(u)}(dv \times dw)$ for all $u > u_{\pi_1}, u \in C_{\pi_1}$ and $\phi \in \Phi$, where $\Pi_{2,3}^{(u)}(A)$ is a measure on the $\sigma$-algebra $B_R$ for every $u > u_{\pi_1}, u \in C_{\pi_1}$.

Obviously $\Pi_{2,3}^{(u)}(A)$ is monotonic in $u > u_\pi, u \in C_{\pi_1}$, because the preliminiting functions in the left-hand side in Condition C are monotonic in $u$ for nonnegative $\phi \in \Phi$.

(a) Due to this property, there exist $\lim_{u' \in C_{\pi_1}, u < u' \to u} \Pi_{2,3}^{(u')}(A) = \Pi_{2,3}^{(u)}(A)$ for every $A \in B_R$, and $u > u_\pi, u \notin C_{\pi_1}$, and also for $u = u_\pi$ if $u_\pi > -\infty, \pi_1(u_\pi) < \infty$.

(b) The following estimates are valid: $\Pi_{2,3}^{(u_1)}(A) - \Pi_{2,3}^{(u_2)}(A) \leq (\pi_1(u_1) - \pi_1(u_2)) \wedge \Pi_{2,3}(A)$; in particular, $\Pi_{2,3}^{(u_1)}(A) \leq \pi_1(u_1) \wedge \Pi_{2,3}(A)$ for any $A \in B_R$, and $u_\pi < u_1 \leq u_2 < \infty$, and also for $u_\pi = u_1 \leq u_2 < \infty$ if $u_\pi > -\infty, \pi_1(u_\pi) < \infty$.

These estimates can be verified by the limiting transition in the corresponding estimates for the preliminiting functions in the left-hand side in Condition C, for nonnegative $\phi \in \Phi$ approximating in a proper way the indicators of sets $A \in B_R$, and $u_\pi < u_1 \leq u_2 < \infty$, $u_1, u_2 \in C_{\pi_1}$. Then by the limiting transition in $u_i' \in C_{\pi_1}, u_i < u_i' \to u_i, i = 1, 2$, these estimates can be obtained for any $A \in B_R$ and $u_\pi < u_1 \leq u_2 < \infty$ and $u_\pi = u_1 \leq u_2 < \infty$ if $u_\pi > -\infty, \pi_1(u_\pi) < \infty$. It follows from these estimates that convergence in statement is uniform with respect to $A \in B_R$. This implies that $\Pi_{2,3}^{(u)}(A)$ is a measure for every $u > u_\pi$ and for $u = u_\pi$ if $u_\pi > -\infty, \pi_1(u_\pi) < \infty$.

(c) These estimates also imply that $\Pi_{2,3}^{(u)}(A)$, as a function in $u$ for every $A \in B_R$, is nondecreasing and right-continuous at any point $u > u_\pi$ and at point $u = u_\pi$ if $u_\pi > -\infty, \pi_1(u_\pi) < \infty$.

For $u > u_\pi$ and for $u = u_\pi$ if $u_\pi > -\infty, \pi_1(u_\pi) < \infty$, we define a measure on the $\sigma$-algebra $B_R$ by the formula

$$(2.6) \qquad \hat\Pi_{2,3}^{(u)}(A) = \Pi_{2,3}(A) - \Pi_{2,3}^{(u)}(A).$$

Let us also define for the corresponding projections $\Pi_2^{(u)}(B) = \Pi_{2,3}^{(u)}(B \times [0,\infty))$ and $\hat\Pi_2^{(u)}(B) = \hat\Pi_{2,3}^{(u)}(B \times [0,\infty))$, which are measures on the Borel $\sigma$-algebra of subsets of $R_1 \setminus \{0\}$, and also define $\Pi_3^{(u)}(C) = \Pi_{2,3}^{(u)}(R_1 \times C)$ and $\hat\Pi_3^{(u)}(C) = \hat\Pi_{2,3}^{(u)}(R_1 \times C)$, which are measures on the Borel $\sigma$-algebra of subsets of $(0,\infty)$.

To be able to write down the representation of the limiting process, we define for $u > u_\pi$ or $u = u_\pi$ if $u_\pi > -\infty$, $\pi_1(u_\pi) < \infty$ and $t > 0$ the characteristic function

$$(2.7) \quad \begin{aligned} &\varphi_{2,3}^{(u)}(t, y, z) \\ &= \exp\Big\{ t\Big( ia^{(u)}y - \frac{1}{2}b^2 y^2 + id^{(u)}z \\ &\qquad\qquad + \int_R \Big( e^{i(yv+zw)} - 1 - \frac{i(yv+zw)}{1+v^2+w^2} \Big) \hat{\Pi}_{2,3}^{(u)}(dv \times dw) \Big) \Big\}, \end{aligned}$$

where

$$(2.8) \quad a^{(u)} = a - \fint_{R_1} \frac{s}{1+s^2} \Pi_2^{(u)}(ds), \qquad d^{(u)} = d - \int_{(0,\infty)} \frac{s}{1+s^2} \Pi_3^{(u)}(ds).$$

Note also that the constants $a, b, d, a^{(u)}$ and $d^{(u)}$ and the measures $\Pi_{2,3}(A)$ and $\Pi_{2,3}^{(u)}(A)$ in (2.6)–(2.8) are determined by Conditions A–C.

It follows from statement (b) that $a^{(u)}, d^{(u)}$ and $\varphi_{2,3}^{(u)}(t, y, z)$ are right-continuous functions in $u > u_\pi$ and for $u = u_\pi$ if $u_\pi > -\infty$, $\pi_1(u_\pi) < \infty$.

Let us also define $\varphi_{2,3}^{(u)}(t, y, z) = 1$ for $u < u_\pi$ or for $u = u_\pi$ if $u_\pi > -\infty$, $\pi_1(u_\pi) = \infty$. Here is a first key result.

THEOREM 2.1. *Let Conditions A–C hold. Then*

$$(2.9) \quad \begin{aligned} &\{(\xi_\varepsilon(t), \gamma_\varepsilon(t), \kappa_\varepsilon(t)), t > 0\} \\ &\Longrightarrow \{(\xi_0(t), \gamma_0(t), \kappa_0(t)), t > 0\} \qquad \text{as } \varepsilon \to 0, \end{aligned}$$

*where $\{(\xi_0(t), \gamma_0(t), \kappa_0(t)), t > 0\}$ is a homogeneous, stochastically continuous Markov process whose trajectories belong to the space $D_0 \times D \times D_{++}$ with probability 1 and transition probabilities that have the hybrid characteristic-distribution form*

$$(2.10) \quad \begin{aligned} &E\{\exp\{i(y(\gamma_0(t+s) - \gamma_0(s)) + z(\kappa_0(t+s) - \kappa_0(s)))\} \\ &\quad \times \chi(\xi_0(t+s) \le u) \,|\, \xi_0(s) = u', \gamma_0(s) = v', \kappa_0(s) = w'\} \\ &= \chi(u' \le u) \exp(-t\pi_1(u)) \cdot \varphi_{2,3}^{(u)}(t, y, z). \end{aligned}$$

PROOF. The method that we use is based on application of the classical central criterion of convergence to distributions of sum processes conditioned in a special way with respect to the corresponding max components. This method was proposed in Chow and Teugels (1979), where asymptotics of joint distributions of maxima and sums of i.i.d. random variables ($\xi_{\varepsilon,k} = \gamma_{\varepsilon,k}, k = 1, 2, \dots$) was investigated for the case of the scale location model in the situation when the random variables belong to the domain of attraction of a stable law. Here, we deal with nonidentical random variables $\xi_{\varepsilon,k}$ and $\gamma_{\varepsilon,k}$ and a general triangular array model. This complicates the consideration. Nevertheless, the method still is the most effective one. It also yields explicit expressions for the corresponding limiting characteristics that are not easy to guess in advance.

By the definition of the processes $\xi_\varepsilon(t)$, $\gamma_\varepsilon(t)$ and $\kappa_\varepsilon(t)$ for any $0 = t_0 < t_1 < \cdots < t_m < \infty$, $-\infty < u_1 < \cdots < u_m < \infty$, $\bar{y} = (y_1, \ldots, y_m)$, $\bar{z} = (z_1, \ldots, z_m) \in R_m$, $m \geq 1$, and for $\varepsilon$ such that $t_1 n_\varepsilon \geq 1$,

$$
\begin{aligned}
(2.11) \quad & E \exp \left\{ i \sum_{l=1}^m (y_l \gamma_\varepsilon(t_l) + z_l \kappa_\varepsilon(t_l)) \right\} \chi(\xi_\varepsilon(t_l) \leq u_l, l = 1, \ldots, m) \\
& = \prod_{l=1}^m (E \exp\{i(y_{l,m}\gamma_{\varepsilon,1} + z_{l,m}\kappa_{\varepsilon,1})\} \chi(\xi_{\varepsilon,1} \leq u_l))^{[t_l n_\varepsilon] - [t_{l-1} n_\varepsilon]},
\end{aligned}
$$

where $y_{l,m} = y_l + \cdots + y_m, z_{l,m} = z_l + \cdots + z_m, l = 1, \ldots, m$.

It follows from (2.11) that (2.9) will result if we can show that when $y, z \in R_1$, for every $u \in C_{\pi_1}$,

$$
\begin{aligned}
(2.12) \quad & (E \exp\{i(y\gamma_{\varepsilon,1} + z\kappa_{\varepsilon,1})\} \chi(\xi_{\varepsilon,1} \leq u))^{n_\varepsilon} \\
& \to \exp(-\pi_1(u)) \varphi_{2,3}^{(u)}(1, y, z) \qquad \text{as } \varepsilon \to 0.
\end{aligned}
$$

This relationship is obvious for the case $u < u_\pi$, since in this case the expression on the right-hand side in (2.12) tends to zero due to Condition A and the expression on the left-hand side in (2.12) is also equal to zero due to the same condition.

If $u_\pi$ is a point of continuity of the function $\pi_1(u)$, then $\pi_1(u_\pi) = \infty$. In this case again the expression on the right-hand side in (2.12) taken for $u = u_\pi$ tends to zero due to Condition A and the expression on the left-hand side in (2.12) taken for $u = u_\pi$ is also equal to zero. So, the only case that needs to be considered is when $u > u_\pi, u \in C_{\pi_1}$.

Obviously

$$
\begin{aligned}
(2.13) \quad & (E \exp\{i(y\gamma_{\varepsilon,1} + z\kappa_{\varepsilon,1})\} \chi(\xi_{\varepsilon,1} \leq u))^{n_\varepsilon} \\
& = (P\{\xi_{\varepsilon,1} \leq u\})^{n_\varepsilon} (E\{\exp\{i(y\gamma_{\varepsilon,1} + z\kappa_{\varepsilon,1})\} \mid \xi_{\varepsilon,1} \leq u\})^{n_\varepsilon}.
\end{aligned}
$$

By Condition A for $u > u_\pi, u \in C_{\pi_1}$,

$$
(2.14) \qquad (P\{\xi_{\varepsilon,1} \leq u\})^{n_\varepsilon} \to \exp(-\pi_1(u)) \qquad \text{as } \varepsilon \to 0.
$$

From (2.13) and (2.14), relation (2.12) will be proved if we show that when $y, z \in R_1$, for every $u > u_\pi, u \in C_{\pi_1}$,

$$
(2.15) \quad (E\{\exp\{i(y\gamma_{\varepsilon,1} + z\kappa_{\varepsilon,1})\} \mid \xi_{\varepsilon,1} \leq u\})^{n_\varepsilon} \to \varphi_{2,3}^{(u)}(1, y, z) \qquad \text{as } \varepsilon \to 0.
$$

For every $\varepsilon > 0$ and $u > u_\pi, u \in C_{\pi_1}$, define sequences of i.i.d. random vectors $\{(\gamma_{\varepsilon,n}^{(u)}, \kappa_{\varepsilon,n}^{(u)}), n = 1, 2, \ldots\}$ such that for $v \in R_1, w \geq 0$,

$$
(2.16) \qquad P\{\gamma_{\varepsilon,1}^{(u)} \leq v, \kappa_{\varepsilon,1}^{(u)} \leq w\} = P\{\gamma_{\varepsilon,1} \leq v, \kappa_{\varepsilon,1} \leq w \mid \xi_{\varepsilon,1} \leq u\}.
$$

With these sequences we can associate their natural sum processes defined by

$$(2.17) \qquad \gamma_\varepsilon^{(u)}(t) = \sum_{k \le tn_\varepsilon} \gamma_{\varepsilon,k}^{(u)}, \qquad \kappa_\varepsilon^{(u)}(t) = \sum_{k \le tn_\varepsilon} \kappa_{\varepsilon,k}^{(u)}, \qquad t \ge 0.$$

For given $u > u_\pi, u \in C_{\pi_1}$, relationship (2.15) is actually equivalent to

$$(2.18) \qquad \begin{aligned} &\{(\gamma_\varepsilon^{(u)}(t), \kappa_\varepsilon^{(u)}(t)), t \ge 0\} \\ &\quad \Longrightarrow \quad \{(\gamma_0^{(u)}(t), \kappa_0^{(u)}(t)), t \ge 0\} \qquad \text{as } \varepsilon \to 0, \end{aligned}$$

where $\{(\gamma_0^{(u)}(t), \kappa_0^{(u)}(t)), t \ge 0\}$ is a homogeneous process with independent increments with the characteristic function $\varphi_{2,3}^{(u)}(t, y, z)$.

It was pointed out in Section 2.1 that Condition B is necessary and sufficient for (2.18) to hold. Of course, all of these conditions should be checked for the random vectors $(\gamma_{\varepsilon,1}^{(u)}, \kappa_{\varepsilon,1}^{(u)})$ rather than for the random vectors $(\gamma_{\varepsilon,1}, \kappa_{\varepsilon,1})$. These conditions should be checked for every $u > u_\pi, u \in C_{\pi_1}$. Comparison of (2.5) and (2.7) shows that we need constants $a^{(u)}, b, d^{(u)}$ and measures $\hat{\Pi}_{2,3}^{(u)}(A)$ to replace constants $a, b, d$ and measures $\Pi_{2,3}(A)$ in these conditions. We deal with all of them in separate steps.

(i) Let us first treat the asymptotic relationship in Condition B(a). Note first that Condition A implies that for every $u > u_\pi$,

$$(2.19) \qquad P\{\xi_{\varepsilon,1} \le u\} \to 1 \qquad \text{as } \varepsilon \to 0.$$

Using Conditions B(a) and C and (2.19) we have, for every $u > u_\pi, u \in C_{\pi_1}$ and every function $\phi \in \Phi$,

$$(2.20) \qquad \begin{aligned} &n_\varepsilon E\phi(\gamma_{\varepsilon,1}^{(u)}, \kappa_{\varepsilon,1}^{(u)}) \\ &= n_\varepsilon \frac{E\chi(\xi_{\varepsilon,1} \le u)\phi(\gamma_{\varepsilon,1}, \kappa_{\varepsilon,1})}{P\{\xi_{\varepsilon,1} \le u\}} \\ &= n_\varepsilon \frac{E\phi(\gamma_{\varepsilon,1}, \kappa_{\varepsilon,1}) - E\chi(\xi_{\varepsilon,1} > u)\phi(\gamma_{\varepsilon,1}, \kappa_{\varepsilon,1})}{P\{\xi_{\varepsilon,1} \le u\}} \\ &\to \int_R \phi(v,w)\Pi_{2,3}(dv \times dw) - \int_R \phi(v,w)\Pi_{2,3}^{(u)}(dv \times dw) \\ &= \int_R \phi(v,w)\hat{\Pi}_{2,3}^{(u)}(dv \times dw) \qquad \text{as } \varepsilon \to 0. \end{aligned}$$

(ii) We turn to the asymptotic relationships given in Conditions B(b) and (c) which have the same structure. We restrict attention to the more general Condition B(b) because the proof of Condition B(c) is analogous.

Using Conditions A and B(b) we have, for every $u > u_\pi, u \in C_{\pi_1}$ and $0 < v_k \to 0$ as $k \to \infty$,

$$
\begin{aligned}
(2.21) \quad & \limsup_{\varepsilon \to 0} |n_\varepsilon E \gamma_{\varepsilon,1} \chi(\xi_{\varepsilon,1} > u, |\gamma_{\varepsilon,1}| \leq v_k)| \\
& \leq \limsup_{\varepsilon \to 0} n_\varepsilon E |\gamma_{\varepsilon,1}| \chi(\xi_{\varepsilon,1} > u, |\gamma_{\varepsilon,1}| \leq v_k) \\
& \leq \limsup_{\varepsilon \to 0} v_k n_\varepsilon P\{\xi_{\varepsilon,1} > u\} = v_k \pi_1(u) \to 0 \qquad \text{as } k \to \infty.
\end{aligned}
$$

We use (2.20) again together with (2.21) to see that for every $u > u_\pi, u \in C_{\pi_1}$ and $0 < v_k < v, v_k \to 0$ such that $\Pi_2^{(u)}(\{\pm v\}) = \Pi_2^{(u)}(\{\pm v_k\}) = 0$,

$$
\begin{aligned}
(2.22) \quad & \lim_{\varepsilon \to 0} n_\varepsilon E \gamma_{\varepsilon,1} \chi(\xi_{\varepsilon,1} > u, |\gamma_{\varepsilon,1}| \leq v) \\
& = \lim_{k \to \infty} \lim_{\varepsilon \to 0} n_\varepsilon E \gamma_{\varepsilon,1} \chi(\xi_{\varepsilon,1} > u, v_k \leq |\gamma_{\varepsilon,1}| \leq v) \\
& = \lim_{k \to \infty} \int_{v_k \leq |s| \leq v} s \Pi_2^{(u)}(ds) = \fint_{|s| \leq v} s \Pi_2^{(u)}(ds).
\end{aligned}
$$

Combining (2.22) and Condition B(b) we get for $u > u_\pi, u \in C_{\pi_1}$ and $v > 0$ such that $\Pi_2^{(u)}(\{\pm v\}) = 0$,

$$
\begin{aligned}
(2.23) \quad & n_\varepsilon E \gamma_{\varepsilon,1}^{(u)} \chi(|\gamma_{\varepsilon,1}^{(u)}| \leq v) \\
& = n_\varepsilon \frac{E \gamma_{\varepsilon,1} \chi(\xi_{\varepsilon,1} \leq u, |\gamma_{\varepsilon,1}| \leq v)}{P\{\xi_{\varepsilon,1} \leq u\}} \\
& = n_\varepsilon \frac{E \gamma_{\varepsilon,1} \chi(|\gamma_{\varepsilon,1}| \leq v) - E \gamma_{\varepsilon,1} \chi(\xi_{\varepsilon,1} > u, |\gamma_{\varepsilon,1}| \leq v)}{P\{\xi_{\varepsilon,1} \leq u\}} \\
& \to a^{(u)}(v) \\
& = a(v) - \fint_{|s| \leq v} s \Pi_{1,2}^{(u)}(ds) \qquad \text{as } \varepsilon \to 0.
\end{aligned}
$$

This relationship enables us to calculate the corresponding constant $a^{(u)}$ in (2.7) that replaces $a$. Indeed, according to relationships (2.20) and (2.23) and the defining formula for the constant $a$, we have

$$
\begin{aligned}
(2.24) \quad a^{(u)} &= a^{(u)}(v) - \fint_{|s| < v} \frac{s^3}{1+s^2} \hat{\Pi}_2^{(u)}(ds) + \int_{|s| > v} \frac{s}{1+s^2} \hat{\Pi}_2^{(u)}(ds) \\
&= a(v) - \fint_{|s| < v} s \Pi_{1,2}^{(u)}(ds) - \fint_{|s| < v} \frac{s^3}{1+s^2} [\Pi_2(ds) - \Pi_2^{(u)}(ds)] \\
&\quad + \int_{|s| > v} \frac{s}{1+s^2} [\Pi_2(ds) - \Pi_2^{(u)}(ds)] \\
&= a - \fint_{|s| < v} s \Pi_2^{(u)}(ds) \\
&\quad + \fint_{|s| < v} \frac{s^3}{1+s^2} \Pi_2^{(u)}(ds) - \int_{|s| > v} \frac{s}{1+s^2} \Pi_2^{(u)}(ds) \\
&= a - \fint_{R_1} \frac{s}{1+s^2} \Pi_2^{(u)}(ds).
\end{aligned}
$$

(iii) Finally, we must check Condition B(d) for the random variables $\gamma_{\varepsilon,1}^{(u)}$.

Note that (2.23) implies in an obvious way that, for $u > u_\pi, u \in C_{\pi_1}$ and $v > 0$,

$$(2.25) \qquad \limsup_{\varepsilon \to 0} n_\varepsilon (E\gamma_{\varepsilon,1}^{(u)} \chi(|\gamma_{\varepsilon,1}^{(u)}| \leq v))^2 = 0.$$

Then, using Conditions A and B(d) we get for $u > u_\pi, u \in C_{\pi_1}$,

$$
\begin{aligned}
(2.26) \quad & \limsup_{\varepsilon \to 0} n_\varepsilon E\gamma_{\varepsilon,1}^2 \chi(\xi_{\varepsilon,1} > u, |\gamma_{\varepsilon,1}| \leq v) \\
& \leq \limsup_{\varepsilon \to 0} \sqrt{n_\varepsilon P\{\xi_{\varepsilon,1} > u\}} \cdot \sqrt{n_\varepsilon E\gamma_{\varepsilon,1}^4 \chi(|\gamma_{\varepsilon,1}| \leq v)} \\
& \leq \limsup_{\varepsilon \to 0} \sqrt{n_\varepsilon P\{\xi_{\varepsilon,1} > u\}} \cdot \limsup_{\varepsilon \to 0} \sqrt{v^2 n_\varepsilon E\gamma_{\varepsilon,1}^2 \chi(|\gamma_{\varepsilon,1}| \leq v)} \\
& \leq \sqrt{\pi_1(u)} \cdot \sqrt{v^2 \limsup_{\varepsilon \to 0} n_\varepsilon E\gamma_{\varepsilon,1}^2 \chi(|\gamma_{\varepsilon,1}| \leq v)} \to 0 \qquad \text{as } 0 < v \to 0.
\end{aligned}
$$

Using (2.19), (2.25), (2.26) and Conditions B(b) and (d) we get, for $u > u_\pi, u \in C_{\pi_1}$,

$$
\begin{aligned}
(2.27) \quad & \lim_{0 < v \to 0} \overline{\overline{\lim}}_{\varepsilon \to 0} n_\varepsilon \operatorname{Var} \gamma_{\varepsilon,1}^{(u)} \chi(|\gamma_{\varepsilon,1}^{(u)}| \leq v) \\
& = \lim_{0 < v \to 0} \overline{\overline{\lim}}_{\varepsilon \to 0} n_\varepsilon E(\gamma_{\varepsilon,1}^{(u)})^2 \chi(|\gamma_{\varepsilon,1}^{(u)}| \leq v) \\
& = \lim_{0 < v \to 0} \overline{\overline{\lim}}_{\varepsilon \to 0} n_\varepsilon \frac{E\gamma_{\varepsilon,1}^2 \chi(|\gamma_{\varepsilon,1}| \leq v) - E\gamma_{\varepsilon,1}^2 \chi(\xi_{\varepsilon,1} > u, |\gamma_{\varepsilon,1}| \leq v)}{P\{\xi_{\varepsilon,1} \leq u\}} \\
& = \lim_{0 < v \to 0} \overline{\overline{\lim}}_{\varepsilon \to 0} n_\varepsilon E\gamma_{\varepsilon,1}^2 \chi(|\gamma_{\varepsilon,1}| \leq v) \\
& = \lim_{0 < v \to 0} \overline{\overline{\lim}}_{\varepsilon \to 0} n_\varepsilon \operatorname{Var} \gamma_{\varepsilon,1}^2 \chi(|\gamma_{\varepsilon,1}| \leq v) = b^2.
\end{aligned}
$$

Note that the constant $b$ does not depend on $u > u_\pi$. By combining the above determinations, the proof is complete. $\quad\square$

### 2.4. *J-convergence of mixed max–sum processes.*

We turn to convergence in $J$-topology. Let $D_k$ be the space of functions on $(0, \infty)$, taking values in $R_k$, that are right-continuous and have no discontinuities of the second kind. The symbol $\{\zeta_\varepsilon(t), t > 0\} \xrightarrow{J} \{\zeta_0(t), t > 0\}$ as $\varepsilon \to 0$ is used to indicate that the processes $\zeta_\varepsilon(t)$, whose trajectories belong to the space $D_k$ with probability 1, converge in Skorokhod $J$-topology to a process $\zeta_0(t)$ on any interval $[t', t'']$, where $0 < t' < t'' < \infty$ are points of stochastic continuity of the process $\{\zeta_0(t), t > 0\}$.

We refer to the books by Billingsley (1968) and Gikhman and Skorokhod (1971) as well as to the articles by Stone (1963) and Lindvall (1973), where one can find the basic definitions and general facts concerning $J$-convergence for random processes on finite and infinite intervals.

Our interest lies in the process $\{\zeta_\varepsilon(t) = (\xi_\varepsilon(t), \gamma_\varepsilon(t), \kappa_\varepsilon(t)), t > 0\}$ which has phase space $R_1 \times R_1 \times [0, \infty)$ and trajectories that by definition belong to the space $D_3$ with probability 1. It is a Markov process. We denote the transition probabilities of this process by $P_\varepsilon((u, v, w), t, t + s, A)$.

The following theorem is our second main result.

THEOREM 2.2. *Let the Conditions* A–C *hold. Then*

$$(2.28) \qquad \{\zeta_\varepsilon(t), t > 0\} \xrightarrow{J} \{\zeta_0(t), t > 0\} \qquad as\ \varepsilon \to 0.$$

PROOF. Let $\{x(t), t > 0\}$ be a function from the space $D_k$ and $0 < T < T' < \infty, c > 0$. Denote the *modulus of compactness for topology J* by

$$\Delta_J(x(\cdot), c, T, T') = \sup_{T \vee (t-c) \leq t' \leq t \leq t'' \leq (t+c) \wedge T'} \min(|x(t') - x(t)|, |x(t'') - x(t)|).$$

Whereas the weak convergence of the processes $\{\zeta_\varepsilon(t), t > 0\}$ has been proven in Theorem 2.1, Theorem 2.2 will follow if we can show that, for all $0 < T < T' < \infty$ and $\delta > 0$,

$$(2.29) \qquad \lim_{c \to 0} \limsup_{\varepsilon \to 0} P\{\Delta_J(\zeta_\varepsilon(\cdot), c, T, T') \geq \delta\} = 0.$$

Note that the first component $\{\xi_\varepsilon(t), t > 0\}$ is a nondecreasing process with probability 1. We use this property to reduce the phase space of the first component to the interval $[h, \infty)$. This is an essential part in the proof of the $J$-compactness relationship (2.29).

Let us choose:

(d) $h > -\infty$ to be a point of continuity of the function $\pi_1(u)$ if $u_\pi = -\infty$;
(e) $h = u_\pi$ if $u_\pi > -\infty$.

We introduce the truncated random variables $\hat{\xi}_{\varepsilon,k}^{(h)} = \xi_{\varepsilon,k} \vee h, k = 1, 2, \ldots,$ and the corresponding max processes

$$(2.30) \qquad \hat{\xi}_\varepsilon^{(h)}(t) = \max_{k \leq 1 \vee tn_\varepsilon} \hat{\xi}_{\varepsilon,k}^{(h)} = \xi_\varepsilon(t) \vee h, \qquad t \geq 0.$$

The three-variate process $\{\hat{\zeta}_\varepsilon^{(h)}(t) = (\hat{\xi}_\varepsilon^{(h)}(t), \gamma_\varepsilon(t), \kappa_\varepsilon(t)), t > 0\}$ has the phase space $[h, \infty) \times R_1 \times [0, \infty)$ and its trajectories belong to the space $D_3$ with probability 1. It is a Markov process which has for $(u, v, w) \in [h, \infty) \times R_1 \times [0, \infty)$ the same transition probabilities $P_\varepsilon((u, v, w), t, t + s, A)$ as the process $\{\zeta_\varepsilon(t), t > 0\}$.

Note that Theorem 2.1 can be applied to max–sum processes $\{\hat{\zeta}_\varepsilon^{(h)}(t), t > 0\}$ because all conditions of Theorem 2.1 are satisfied. The only difference is that in the current case the corresponding limiting functions and measures in Conditions A–C should be changed. We introduce new functions indexed by a lower index $h$ as follows: $\pi_{1;h}(u) = \pi_1(u)$, $\Pi_{2,3;h}(A) = \Pi_{2,3}(A)$ for $u \geq h$, while $\pi_{1;h}(u) = \infty$, $\Pi_{2,3;h}(A) = \Pi_{2,3}(A)$ for $u < h$. Note that only in the case when $u_\pi = -\infty$ are the changes genuine, while in the case $u_\pi > -\infty$ and $h = u_\pi$ the new functions coincide with the old ones.

According to Theorem 2.1, the following relationship holds:

$$(2.31) \qquad \{\hat{\zeta}_\varepsilon^{(h)}(t), t > 0\} \quad \Longrightarrow \quad \{\hat{\zeta}_0^{(h)}(t), t \geq 0\} \qquad \text{as } \varepsilon \to 0.$$

This relationship also follows directly from (2.9) since according to (2.30) the random variable $\hat{\xi}_\varepsilon^{(h)}(t)$ is a continuous function of the random variable $\xi_\varepsilon(t)$ for every $t > 0$ and $h \in R_1$.

The limiting process $\{\hat{\zeta}_0^{(h)}(t), t > 0\}$ is fully similar to the process described in Theorem 2.1 with the only changes that $h$ is the lower index in the limiting characteristics. Moreover, it can be easily shown that the process $\{\hat{\zeta}_0^{(h)}(t), t > 0\}$ can be constructed from the process $\{\zeta_0(t), t > 0\}$ by simple truncation of the first component of this process, that is, $\{\hat{\zeta}_0^{(h)}(t) = (\hat{\xi}_0^{(h)}(t), \gamma_0(t), \kappa_0(t)), t > 0\}$, where $\hat{\xi}_0^{(h)}(t) = \xi_0(t) \vee h, t > 0$. Note again that in the case $u_\pi = -\infty$ the truncation is genuine while in the cases $u_\pi > -\infty$ and $h = u_\pi$ process $\hat{\xi}_0^{(h)}(t) \equiv \xi_0(t), t > 0$.

Let us now use the inequality

$$(2.32) \qquad \begin{aligned} &P\{\Delta_J(\zeta_\varepsilon(\cdot), c, T, T') \geq 2\delta\} \\ &\qquad \leq P\{\Delta_J(\hat{\zeta}_\varepsilon^{(h)}(\cdot), c, T, T') \geq \delta\} + P\Big\{\sup_{t \in [T,T']} |\hat{\xi}_\varepsilon^{(h)}(t) - \xi_\varepsilon(t)| \geq \delta\Big\}. \end{aligned}$$

Obviously,

$$(2.33) \qquad P\Big\{\sup_{t \in [T,T']} |\hat{\xi}_\varepsilon^{(h)}(t) - \xi_\varepsilon(t)| \geq \delta\Big\} \leq P\{\xi_\varepsilon(T) \leq h - \delta\}.$$

One can always choose $\delta/2 \leq \delta_h \leq \delta$ in such a way that the point $h - \delta_h$ is a point of continuity of function $\pi_1(u)$. Then we get, for $\varepsilon$ such that $n_\varepsilon T \geq 1$,

$$(2.34) \qquad \begin{aligned} &\lim_{h \to -\infty} \limsup_{\varepsilon \to 0} P\{\xi_\varepsilon(T) \leq h - \delta\} \\ &\qquad \leq \lim_{h \to -\infty} \limsup_{\varepsilon \to 0} P\{\xi_\varepsilon(T) \leq h - \delta_h\} \\ &\qquad = \lim_{h \to -\infty} \limsup_{\varepsilon \to 0} P\{\xi_{\varepsilon,1} \leq h - \delta_h\}^{[n_\varepsilon T]} \\ &\qquad = \lim_{h \to -\infty} e^{-\pi_1(h - \delta_h)T} = 0. \end{aligned}$$

Inequalities (2.33) and relationship (2.34) imply that for any $0 < T < T' < \infty$ and $\delta > 0$,

$$(2.35) \qquad \lim_{h \to -\infty} \limsup_{\varepsilon \to 0} P\Big\{\sup_{t \in [T,T']} |\hat{\xi}_\varepsilon^{(h)}(t) - \xi_\varepsilon(t)| \geq \delta\Big\} = 0.$$

In the case $u_\pi > -\infty$, the internal limiting expression on the left-hand side in (2.34) is equal to zero. In this case additional external limit transitions given in (2.34) and (2.35) are not required.

Relationships (2.32) and (2.35) imply that (2.29) will follow if we show that for any $h$, chosen as described above, and any $0 < T < T' < \infty$ and $\delta > 0$,

$$(2.36) \qquad \lim_{c \to 0} \limsup_{\varepsilon \to 0} P\{\Delta_J(\hat{\zeta}_\varepsilon^{(h)}(\cdot), c, T, T') \geq \delta\} = 0.$$

Let us define

$$\alpha_\varepsilon(h, c, T, T', \delta)$$
$$= \sup_{u \geq h, -\infty < v < \infty, w \geq 0} \sup_{T \leq t \leq t + s \leq t + c \leq T'} P_\varepsilon((u, v, w), t, t + s, S_\delta((u, v, w))),$$

where $S_\delta((u, v, w)) = \{(u', v', w') : (|u - u'|^2 + |v - v'|^2 + |w - w'|^2)^{1/2} > \delta\}$.

We showed in (2.31) weak convergence of processes $\{\hat{\zeta}_\varepsilon^{(h)}(t), t > 0\}$. As is known [see Gikhman and Skorokhod (1971)], in this case (2.36) follows from the following relationship that should be proved for any $0 < T < T' < \infty$ and $\delta > 0$:

$$(2.37) \qquad \lim_{c \to 0} \limsup_{\varepsilon \to 0} \alpha_\varepsilon(h, c, T, T', \delta) = 0.$$

We now exploit the fact that the process $\{\xi_\varepsilon(t), t > 0\}$ is nondecreasing and that the two processes $\{\gamma_\varepsilon(t), t > 0\}$ and $\{\kappa_\varepsilon(t), t > 0\}$ are processes with independent increments. We get the estimate

$$(2.38) \quad \begin{aligned} &\alpha_\varepsilon(h, c, T, T', 3\delta) \\ &\leq \sup_{u \geq h} \sup_{T \leq t \leq t + s \leq t + c \leq T'} P\{\xi_\varepsilon(t + s) - \xi_\varepsilon(t) > \delta | \xi_\varepsilon(t) = u\} \\ &\quad + \sup_{T \leq t \leq t + s \leq t + c \leq T'} (P\{|\gamma_\varepsilon(t + s) - \gamma_\varepsilon(t)| > \delta\} \\ &\qquad\qquad\qquad\qquad + P\{\kappa_\varepsilon(t + s) - \kappa_\varepsilon(t) > \delta\}) \\ &\leq \sup_{u \geq h} \sup_{T \leq t \leq t + s \leq t + c \leq T'} (1 - P\{\xi_{\varepsilon,1} \leq u + \delta\}^{[n_\varepsilon(t + s)] - [n_\varepsilon t]}) \\ &\quad + \sup_{T \leq t \leq t + s \leq t + c \leq T'} ([n_\varepsilon(t + s)] - [n_\varepsilon t]) \\ &\qquad \times (P\{|\gamma_{\varepsilon,1}| > \delta\} + |E\gamma_{\varepsilon,1}\chi(|\gamma_{\varepsilon,1}| \leq \delta)| \\ &\qquad\quad + \operatorname{Var}\gamma_{\varepsilon,1}\chi(|\gamma_{\varepsilon,1}| \leq \delta) + P\{\kappa_{\varepsilon,1} > \delta\} + E\kappa_{\varepsilon,1}\chi(\kappa_{\varepsilon,1} \leq \delta)) \\ &\leq 1 - (P\{\xi_{\varepsilon,1} \leq h + \delta\})^{n_\varepsilon c} \\ &\quad + cn_\varepsilon(P\{|\gamma_{\varepsilon,1}| > \delta\} + |E\gamma_{\varepsilon,1}\chi(|\gamma_{\varepsilon,1}| \leq \delta)| \\ &\qquad\quad + \operatorname{Var}\gamma_{\varepsilon,1}\chi(|\gamma_{\varepsilon,1}| \leq \delta) + P\{\kappa_{\varepsilon,1} > \delta\} + E\kappa_{\varepsilon,1}\chi(\kappa_{\varepsilon,1} \leq \delta)). \end{aligned}$$

We are now in a position to use the truncation of the phase space described above, and Conditions A–C.

Indeed, for every $h$, chosen according to choice (d) or (e), $\delta/2 \le \delta_h \le \delta$ can be chosen such that

$$(2.39) \qquad \pi_1(h + \delta_h) < \infty.$$

Conditions A–C and (2.39) applied to (2.38) yield

$$
\begin{aligned}
(2.40) \quad & \lim_{c \to 0} \limsup_{\varepsilon \to 0} \alpha_\varepsilon(h, c, T, T', 3\delta) \\
& \le \lim_{c \to 0} \limsup_{\varepsilon \to 0} (1 - (P\{\xi_{\varepsilon,1} \le h + \delta_h\})^{n_\varepsilon c}) \\
& \quad + \lim_{c \to 0} c \cdot \limsup_{\varepsilon \to 0} n_\varepsilon (P\{|\gamma_{\varepsilon,1}| > \delta\} + |E\gamma_{\varepsilon,1}\chi(|\gamma_{\varepsilon,1}| \le \delta)| \\
& \qquad\qquad\qquad + \operatorname{Var}\gamma_{\varepsilon,1}\chi(|\gamma_{\varepsilon,1}| \le \delta) \\
& \qquad\qquad\qquad + P\{\kappa_{\varepsilon,1} > \delta\} + E\kappa_{\varepsilon,1}\chi(\kappa_{\varepsilon,1} \le \delta)) \\
& = \lim_{c \to 0} (1 - \exp(-\pi_1(h + \delta_h)c)) = 0.
\end{aligned}
$$

The proof is complete. $\quad \square$

## 3. Mixed max–sum processes with renewal stopping.
In this section we deal with general conditions for the weak and functional convergence of the mixed max–sum processes with renewal stopping.

3.1. *Joint weak convergence of max–sum processes and renewal stopping processes.* To start the discussion we need some results about the renewal counting processes. We introduce the renewal stopping processes $\{\tau_\varepsilon(t) = \sup(s : \kappa_\varepsilon(s) \le t), t > 0\}$, where $\kappa_\varepsilon(t) = \sum_{k \le t n_\varepsilon} \kappa_{\varepsilon,k}, t \ge 0$, have been introduced in Section 1. We interpret the sequence $\kappa_{\varepsilon,k}$ as the times between renewals, and hence $\tau_\varepsilon(t)n_\varepsilon - 1$ is the number of renewals in the interval $[0, t n_\varepsilon]$.

We assume the basic Condition B. To exclude the trivial case where the process $\{\kappa_0(t) = 0, t > 0\}$, we also assume the following condition on the measure $\Pi_3(C)$ and the quantity $c$ from Condition B:

CONDITION D. We have $c > 0$ or $\Pi_3((v, \infty)) > 0$ for some $v > 0$.

If we look at the inverse process $\{\tau_0(t) = \sup(s : \kappa_0(s) \le t), t > 0\}$, then Condition D implies that $\kappa_0(t) \xrightarrow{P} \infty$ as $t \to \infty$ and therefore $\tau_0(t)$ is an a.s. finite random variable for every $t > 0$. Conditions B and D obviously imply that, for every $\varepsilon$ small enough (say $\varepsilon \le \varepsilon_0$), $\kappa_\varepsilon(t) \xrightarrow{P} \infty$ as $t \to \infty$ and therefore also the renewal stopping process $\tau_\varepsilon(t)$ is an a.s. finite random variable for every $t > 0$. Furthermore, by definition, the trajectories of the processes $\{\tau_\varepsilon(t), t > 0\}$ a.s. belong to the space $D_{++}$ for every $\varepsilon \le \varepsilon_0$.

Let us denote by $V$ the set of points of stochastic continuity of the process $\{\tau_0(t), t > 0\}$. This process is stochastically continuous, that is, $V = (0, \infty)$ if one of the following conditions hold: (i) $c > 0$, (ii) $\Pi_3((v, \infty)) \to \infty$ as $v \to 0$ or (iii) $\Pi_3((v, \infty))$ is a continuous function. If all three conditions (i)–(iii) are violated, then the set $V$ is $(0, \infty)$ excluding perhaps some countable or finite set. Namely, the set $\overline{V} = (0, \infty) \setminus V = \{v_1 l_1 + \cdots + v_m l_m : l_1, \ldots, l_m = 0, 1, \ldots, l_1 + \cdots + l_m \geq 1, m \geq 1\}$, where $\{v_1, v_2, \ldots\}$ is the set of discontinuity points of the function $\Pi_3((v, \infty))$. The process $\{\tau_0(t), t > 0\}$ is a.s. continuous if the process $\{\kappa_0(t), t > 0\}$ is strictly monotone, that is, if at least one of the conditions (i) or (ii) holds.

Preparing the results for the mixed max–sum processes with renewal stopping, we formulate conditions for the joint weak convergence of max–sum processes and renewal stopping processes.

LEMMA 3.1.  *Let Conditions* A–D *be satisfied. Then*

$$
\begin{align}
(3.1) \qquad & \{(\zeta_\varepsilon(t), \tau_\varepsilon(s)), (t, s) \in (0, \infty) \times V\} \notag \\
& \implies \{(\zeta_0(t), \tau_0(s)), (t, s) \in (0, \infty) \times V\} \qquad \text{as } \varepsilon \to 0. \notag
\end{align}
$$

PROOF.  Use the definition of the processes $\{\zeta_\varepsilon(t), t > 0\}$ and $\{\tau_\varepsilon(t), t > 0\}$ to write that, for any $0 < s_1 < \cdots < s_m$, $0 < t_1 < \cdots < t_m$ and $u_l, v_l, w_l \in R_1, l = 1, 2, \ldots, m, m \geq 1$,

$$
\begin{align}
(3.2) \qquad & P\{\xi_\varepsilon(t_l) \leq u_l, \gamma_\varepsilon(t_l) \leq v_l, \kappa_\varepsilon(t_l) \leq w_l, \tau_\varepsilon(t_l) > s_l, l = 1, \ldots, m\} \notag \\
& = P\{\xi_\varepsilon(t_l) \leq u_l, \gamma_\varepsilon(t_l) \leq v_l, \kappa_\varepsilon(t_l) \leq w_l, \kappa_\varepsilon(s_l) \leq t_l, l = 1, \ldots, m\}. \notag
\end{align}
$$

Choose some countable set of points $X = \{x_1, x_2, \ldots\} \subset (0, \infty)$ dense in $(0, \infty)$. Since any distribution function has at most a countable set of discontinuity points, we can choose a countable set $Y = \{y_1, y_2, \ldots\} \subset V$ dense in $(0, \infty)$ such that $P\{\kappa_0(x_i) = y_j\}$ for all $i, j \geq 1$, and then a countable set $Z = \{z_1, z_2, \ldots\} \subset R_1$ dense in $R_1$ for which $P\{\xi_0(y_j) = z_k\} = P\{\gamma_0(y_j) = z_k\} = P\{\kappa_0(y_j) = z_k\} = 0$ for all $j, k \geq 1$. By Theorem 2.1 and (3.2) we have, for points $s_l \in X, t_l \in Y, u_l, v_l, w_l \in Z, l = 1, \ldots, m, m \geq 1$,

$$
\begin{align}
& P\{\xi_\varepsilon(t_l) \leq u_l, \gamma_\varepsilon(t_l) \leq v_l, \kappa_\varepsilon(t_l) \leq w_l, \tau_\varepsilon(t_l) > s_l, l = 1, \ldots, m\} \notag \\
(3.3) \qquad & \to P\{\xi_0(t_l) \leq u_l, \gamma_0(t_l) \leq u_l, \kappa_0(t_l) \leq w_l, \tau_0(t_l) > s_l, l = 1, \ldots, m\} \notag \\
& \hspace{8cm} \text{as } \varepsilon \to 0. \notag
\end{align}
$$

Note that Condition D implies that all random variables in (3.3) are proper for $\varepsilon$ small enough. Taking into account that the weak convergence of distribution functions of random vectors follows from their convergence on some countable set everywhere dense in the corresponding phase space, we get from (3.3) that

$$
(3.4) \qquad \{(\zeta_\varepsilon(t), \tau_\varepsilon(t)), t \in Y\} \implies \{(\zeta_0(t), \tau_0(t)), t \in Y\} \qquad \text{as } \varepsilon \to 0.
$$

The processes $\{\zeta_\varepsilon(t), t > 0\}$ $J$-converge while the processes $\{\tau_\varepsilon(t), t > 0\}$ are monotonic. Therefore, (3.4) can be extended by obvious arguments to (3.1). $\square$

3.2. *Weak convergence of max–sum processes with renewal stopping.* We are in a position to formulate and prove our third main result.

THEOREM 3.1. *Let Conditions* A–D *hold. Then,*

$$\{\zeta_\varepsilon(\tau_\varepsilon(t)) = (\xi_\varepsilon(\tau_\varepsilon(t)), \gamma_\varepsilon(\tau_\varepsilon(t)), \kappa_\varepsilon(\tau_\varepsilon(t))), t \in V\}$$

(3.5)
$$\implies \quad \{\zeta_0(\tau_0(t)) = (\xi_0(\tau_0(t)), \gamma_0(\tau_0(t)), \kappa_0(\tau_0(t))), t \in V\}$$

$$as\ \varepsilon \to 0.$$

PROOF. We are going to use Theorem 2.7.1 from the book by Silvestrov [(1974), page 82]. This theorem provides the proper conditions for the weak convergence for cadlag processes stopped in Markov type moments. Here, the Markov processes are $\{\hat{\zeta}_\varepsilon^{(h)}(t)\} = \{(\hat{\xi}_\varepsilon^{(h)}(t), \gamma_\varepsilon(t), \kappa_\varepsilon(t)), t > 0\}$ while the stopping moments are given by $\tau_\varepsilon(t), t > 0$. By definition, for every $t > 0$ the random variable $\tau_\varepsilon(t)$ is a Markov moment for the Markov process $\{\hat{\xi}_\varepsilon^{(h)}(t), t > 0\}$. Moreover, also by definition, $P\{\tau_\varepsilon(t) > 0\} = 1, t > 0$, for all $\varepsilon \geq 0$.

As was mentioned above, Theorem 2.1 can be applied to hybrid max–sum processes $\{\hat{\zeta}_\varepsilon^{(h)}(t), t > 0\}$. Also Lemma 3.1 can be applied to these processes and renewal stopping processes $\tau_\varepsilon(t), t > 0$. All conditions of Lemma 3.1 are satisfied. According to Lemma 3.1 the following relationship holds:

(3.6)
$$\{(\hat{\zeta}_\varepsilon^{(h)}(t), \tau_\varepsilon(s)), (t, s) \in (0, \infty) \times V\}$$
$$\implies \{(\hat{\zeta}_0^{(h)}(t), \tau_0(s)), (t, s) \in (0, \infty) \times V\} \quad as\ \varepsilon \to 0.$$

Again it is useful to note that (3.6) also follows directly from (3.1) because $\hat{\xi}_\varepsilon^{(h)}(t) = \xi_\varepsilon(t) \vee h$ is a continuous function of random variable $\xi_\varepsilon(t)$ for every $t > 0$ and $h \in R_1$.

As was mentioned above, the process $\{\hat{\zeta}_0^{(h)}(t), t > 0\}$ can be constructed from the process $\{\zeta_0(t), t > 0\}$ by simple truncation of the first component of this process, that is, $\{\hat{\zeta}_0^{(h)}(t) = (\hat{\xi}_0^{(h)}(t), \gamma_0(t), \kappa_0(t)), t > 0\}$, where $\hat{\xi}_0^{(h)}(t) = \xi_0(t) \vee h, t > 0$. As far as the limiting renewal stopping process is concerned, it can again be defined as $\{\tau_0(t) = \sup(s : \kappa_0(s) \leq t), t > 0\}$.

Next, define

$$\beta_\varepsilon(h, c, t, \delta)$$
$$= \sup_{u \geq h, -\infty < v < \infty, w \geq 0} \sup_{t - c \leq s \leq s + q \leq t + c} P_\varepsilon((u, v, w), s, s + q, S_\delta((u, v, w))).$$

It is easy to see that, for any $0 < T < t < T'$ and for all $c > 0$ such that $c < t, T \leq t - c, t + c \leq T'$,

(3.7)
$$\beta_\varepsilon(h, c, t, \delta) \leq \alpha_\varepsilon(h, 2c, T, T', \delta).$$

Because of this, (2.37) implies that for all $t > 0$,

$$(3.8) \qquad \lim_{c \to 0} \limsup_{\varepsilon \to 0} \beta_\varepsilon(h, c, t, \delta) = 0.$$

According to Theorem 2.7.1 from Silvestrov (1974), (3.6) and (3.8) imply that

$$\{\hat{\zeta}_\varepsilon^{(h)}(\tau_\varepsilon(t)) = (\xi_\varepsilon^{(h)}(\tau_\varepsilon(t)), \gamma_\varepsilon(\tau_\varepsilon(t)), \kappa_\varepsilon(\tau_\varepsilon(t))), t \in V\}$$

$$(3.9) \qquad \Longrightarrow \quad \{\hat{\zeta}_0^{(h)}(\tau_0(t)) = (\hat{\xi}_0^{(h)}(\tau_0(t)), \gamma_0(\tau_0(t)), \kappa_0(\tau_0(t))), t \in V\}$$

$$\text{as } \varepsilon \to 0.$$

We are now prepared to complete the proof of the theorem. Let us use the following inequality, which holds for any $t > 0$ and $0 < T < T' < \infty$:

$$(3.10) \qquad \begin{aligned} & P\{|\hat{\xi}_\varepsilon^{(h)}(\tau_\varepsilon(t)) - \xi_\varepsilon(\tau_\varepsilon(t))| \geq \delta\} \\ & \leq P\{\tau_\varepsilon(t) < T\} + P\{\tau_\varepsilon(t) > T'\} + P\Big\{\sup_{t \in [T, T']} |\hat{\xi}_\varepsilon^{(h)}(t) - \xi_\varepsilon(t)| \geq \delta\Big\}. \end{aligned}$$

For any $t \in V$, points $0 < T < T' < \infty$ can be chosen in such a way that they are continuity points of the distribution functions of random variable $\tau_0(t)$. Moreover, for an arbitrary $\sigma > 0$, the points $T$ and $T'$ can be chosen, respectively, so small and so large that $P\{\tau_0(t) < T\} + P\{\tau_0(t) > T'\} \leq \sigma$.

We use inequality (3.10), relationship (2.35) and the remark concerning the choice of $T, T'$ made above to conclude that for any $0 < t \in V$ and $\delta > 0$,

$$(3.11) \qquad \lim_{h \to -\infty} \limsup_{\varepsilon \to 0} P\{|\hat{\xi}_\varepsilon^{(h)}(\tau_\varepsilon(t)) - \xi_\varepsilon(\tau_\varepsilon(t))| \geq \delta\} = 0.$$

As was already mentioned above, the process $\{\hat{\zeta}_0^{(h)}(t), t > 0\}$ can be constructed from the processes $\{\zeta_0(t), t > 0\}$ by truncation of the first component of this process, that is, $\{\hat{\zeta}_0^{(h)}(t) = (\hat{\xi}_0^{(h)}(t), \gamma_0(t), \kappa_0(t)), t > 0\}$, where $\hat{\xi}_0^{(h)}(t) = \xi_0(t) \vee h, t > 0$, and the limiting renewal stopping process can be defined as $\{\tau_0(t) = \sup(s : \kappa_0(s) \leq t), t > 0\}$. The representation of the limiting processes described above yields that

$$(3.12) \qquad \{\hat{\zeta}_0^{(h)}(\tau_0(t)), t \in V\} \quad \Longrightarrow \quad \{\zeta_0(\tau_0(t)), t \in V\} \qquad \text{as } \varepsilon \to 0.$$

Note that in the case $u_\pi > -\infty$ and $h = u_\pi$, the representation described above yields that $\hat{\xi}_0^{(h)}(t) \equiv \xi_0(t), t > 0$, so that in the sequel $\hat{\xi}_0^{(h)}(\tau_0(t)) \equiv \zeta_0(\tau_0(t)), t > 0$. For this reason, the additional external limit transitions given in (3.11) and (3.12) are not required.

Relationships (3.9), (3.11) and (3.12) imply (3.5).  $\square$

3.3. *J-convergence of renewal processes.* In this section we treat conditions for *J*-convergence of max–sum renewal processes with renewal stopping. We again go through a two-step procedure.

There are two cases when Condition B provides the *J*-convergence relationship

$$(3.13) \qquad \{\tau_\varepsilon(t), t > 0\} \xrightarrow{J} \{\tau_0(t), t > 0\} \qquad \text{as } \varepsilon \to 0.$$

The first case is when the following condition holds:

CONDITION $D_1$. We have $c > 0$ or $\Pi_3((v, \infty)) \to \infty$ as $v \to 0$.

In this case $\{\kappa_0(t), t > 0\}$ is an a.s. strictly monotonic process and therefore $\tau_0(t)$ is an a.s. continuous process. Also, by the definition, the processes $\{\tau_\varepsilon(t), t > 0\}$ and $\{\tau_0(t), t > 0\}$ are nondecreasing.

Since Condition $D_1$ is stronger than Condition D, Lemma 3.1 guarantees that the processes $\tau_\varepsilon(t)$ weakly converge to the processes $\tau_0(t)$ on the set $V$, which is in this case the interval $(0, \infty)$. As is known [see Billingsley (1968)], the weak convergence of monotone processes to a continuous process on the dense set implies their *J*-convergence. Thus (3.13) holds.

The second case is when the following condition holds:

CONDITION $D_2$. As $\varepsilon \to 0$, $n_\varepsilon P\{\kappa_{\varepsilon,1} > 0\} \to \Pi_3((0, \infty))$, where $0 < \Pi_3((0, \infty)) < \infty$.

In this case the process $\kappa_0(t)$ is a compound Poisson process.

Note first that *J*-convergence of the processes $\tau_\varepsilon(t)$ can be derived from conditions for *J*-convergence on monotonic processes as given in Silvestrov (1974) and in Jacod and Shiryaev (1987).

Alternatively, *J*-compactness of the processes $\tau_\varepsilon(t)$ can be obtained by direct estimates for the modulus of *J*-compactness $\Delta_J(\tau_\varepsilon(\cdot), c, T, T')$ obtained by using the fact that the processes $\kappa_\varepsilon(t)$ and $\kappa_0(t)$ are both step processes with a finite number of jumps in any finite interval. Take $0 < T < T' < \infty$ as two points that do not belong to the set $V$. Define $\theta_\varepsilon[T, T']$ to be the minimal length of the intervals between jump points of the process $\tau_\varepsilon(t)$ in the interval $[T, T']$ if there are at least two such points; otherwise put $\theta_\varepsilon[T, T'] = T' - T$. Conditions B and $D_2$ give the joint weak convergence of the consecutive moments and values of the jumps for the processes $\kappa_\varepsilon(t)$ to the corresponding functionals for the processes $\kappa_0(t)$. However, this implies the same joint weak convergence for the consecutive moments and values of the jumps for the inverse processes $\tau_\varepsilon(t)$ and $\tau_0(t)$. It easily follows from the remark made above that the random variables $\theta_\varepsilon(T, T') \Rightarrow \theta_0[T, T']$ as $\varepsilon \to 0$ and that the random variable $\theta_0[T, T'] > 0$ with probability 1. Relationship (3.14) now follows from these two observations and the obvious fact that $\{\Delta_J(\tau_\varepsilon(\cdot), c, T, T') > 0\} \subseteq \{\theta_\varepsilon[T, T'] \leq c\}$. Combining all of this information together tells us that

$$
\begin{aligned}
(3.14) \qquad &\lim_{c \to 0} \limsup_{\varepsilon \to 0} P\{\Delta_J(\tau_\varepsilon(\cdot), c, T, T') \geq \delta\} \\
&\leq \lim_{c \to 0} \limsup_{\varepsilon \to 0} P\{\theta_\varepsilon[T, T'] \leq c\} = 0.
\end{aligned}
$$

It is clear that Condition $D_2$ implies Condition D. Thanks to Lemma 3.1, the processes $\tau_\varepsilon(t)$ weakly converge to the processes $\tau_0(t)$ on the set $V$, which is in this case the interval $(0, \infty)$ except at most a countable set. Therefore, Conditions B and $D_2$ imply the requested $J$-convergence (3.13).

We refer to the recent work by Silvestrov (2002), where one can find a more detailed presentation of the proof. To indicate the relevance of some of the above conditions, an example is given where the process $\kappa_0(t)$ is a compound Poisson process and where the processes $\tau_\varepsilon(t)$ weakly converge but do not $J$-converge when Condition B holds but Condition $D_2$ fails.

### 3.4. $J$-convergence of max–sum processes with renewal stopping.

Everything is ready to prove the fourth key result of the paper.

THEOREM 3.2. *Let the Conditions* A–C *hold together with Conditions* $D_1$ *or* $D_2$. *Then*

$$
(3.15) \qquad \{\zeta_\varepsilon(\tau_\varepsilon(t)), t > 0\} \xrightarrow{J} \{\zeta_0(\tau_0(t)), t > 0\} \qquad as \ \varepsilon \to 0.
$$

PROOF. The weak convergence of the processes $\{\zeta_\varepsilon(\tau_\varepsilon(t)), t \in V\}$ was proved in Theorem 3.1. Recall that the set $V$ is dense in the interval $(0, \infty)$. Therefore, we need only to check $J$-compactness. More precisely, for any $0 < T < T' < \infty$, we need that

$$
(3.16) \qquad \lim_{c \to 0} \limsup_{\varepsilon \to 0} P\{\Delta_J(\zeta_\varepsilon(\tau_\varepsilon(\cdot)), c, T, T') \geq \delta\} = 0.
$$

We go through the proof under the two alternatives.

(i) Assume first that Condition $D_1$ holds. Then we can use results given in Silvestrov (1972b, 1973, 1974), where conditions for $J$-convergence and $J$-compactness of the composition of cadlag processes have been obtained. Actually, we can apply Theorem 2.2.3 from Silvestrov [(1974), page 96] to processes $\{\zeta_\varepsilon(t), t \geq 0\}$ and $\{\tau_\varepsilon(t), t \geq 0\}$. Whereas we need to consider intervals of the form $[T, T']$ rather than $[0, T]$, we repeat here the necessary estimates.

Let $\{x(t), t > 0\}$ be a function from the space $D_1$ and let $0 < T' < T'' < \infty, c > 0$. The *modulus of compactness for the uniform topology* $U$ is defined and denoted by

$$
\Delta_U(x(\cdot), c, T, T') = \sup_{|t' - t''| \leq c, T \leq t' \leq t'' \leq T'} (|x(t') - x(t'')|).
$$

For any $0 < S < S' < \infty$ and $c, C > 0$ we have

$$
\begin{aligned}
(3.17) \quad & P\{\Delta_J(\zeta_\varepsilon(\tau_\varepsilon(\cdot)), c, T, T') \geq \delta\} \\
& \leq P\{\tau_\varepsilon(T) < S\} + P\{\tau_\varepsilon(T') > S'\} \\
& \quad + P\{\Delta_U(\tau_\varepsilon(\cdot), c, T, T') \geq C\} + P\{\Delta_J((\zeta_\varepsilon(\cdot)), C, S, S') \geq \delta\}.
\end{aligned}
$$

Now, $J$-convergence of the processes $\{\zeta_\varepsilon(t), t \geq 0\}$ was proved in Theorem 2.2. As was mentioned above, the weak convergence of the processes $\{\tau_\varepsilon(t), t \in V\}$ to a.s. continuous process $\{\tau_0(t), t \in V\}$ implies the $J$-convergence. Since the limiting process is a.s. continuous, $J$-convergence is equivalent to convergence of these processes in the uniform topology. Also by Theorem 2.1, the processes $\{\tau_\varepsilon(t), t \in V\}$ converge weakly.

It therefore suffices to prove (3.16) only for those $0 < T < T' < \infty$ for which $T, T' \in V$. Choose $S$ and $S'$ continuity points for the distribution functions of the random variables $\tau_\varepsilon(T)$ and $\tau_\varepsilon(T')$, respectively. By using convergence of processes $\{\tau_\varepsilon(t), t \geq 0\}$ in uniform topology, we get, from (3.17),

$$
\begin{aligned}
(3.18) \quad & \lim_{c \to 0} \limsup_{\varepsilon \to 0} P\{\Delta_J(\zeta_\varepsilon(\tau_\varepsilon(\cdot)), c, T, T') \geq \delta\} \\
& \leq P\{\tau_0(T) < S\} + P\{\tau_0(T') > S'\} \\
& \quad + \limsup_{\varepsilon \to 0} P\{\Delta_J((\zeta_\varepsilon(\cdot)), C, S, S') \geq \delta\}.
\end{aligned}
$$

Due to $J$-convergence of processes $\{\zeta_\varepsilon(t), t \geq 0\}$, the expression on the right-hand side can be made less than any $\sigma > 0$ by first choosing $S$ small enough and $S'$ large enough and then by choosing $C$ small enough. This proves (3.16).

(ii) Next assume that Condition $D_2$ holds. In this simpler case, the process $\{\zeta_\varepsilon(\tau_\varepsilon(t)), t > 0\}$ has stepwise trajectories since the internal stopping processes $\{\tau_\varepsilon(t), t > 0\}$ has such trajectories. Moreover, both processes have the same jump points. For this reason $\{\Delta_J(\zeta_\varepsilon(\tau_\varepsilon(\cdot)), c, T, T') > 0\} \subseteq \{\theta_\varepsilon[T, T'] \leq c\}$.

Using this relationship we get as with (3.14) that

$$
\begin{aligned}
(3.19) \quad & \lim_{c \to 0} \limsup_{\varepsilon \to 0} P\{\Delta_J(\zeta_\varepsilon(\tau_\varepsilon(\cdot)), c, T, T') \geq \delta\} \\
& \leq \lim_{c \to 0} \limsup_{\varepsilon \to 0} P\{\theta_\varepsilon[T, T'] \leq c\} = 0.
\end{aligned}
$$

The proof is complete. $\square$

**4. Examples and remarks.** In this section we discuss some variations of results given above and consider examples.

4.1. *Asymptotically independent components.* The max process $\{\xi_0(t), t \geq 0\}$ and the sum process $\{\gamma_0(t), t \geq 0\}$ are independent if and only if the function $\Pi_2^{(u)}((-\infty, -v] \cup (v, \infty)) \equiv 0$ for $u > u_\pi$, $v > 0$. In this case, $a^{(u)} \equiv a$ and the measure $\hat{\Pi}_2^{(u)}(B) \equiv \Pi_2(B)$ for every $u \in R_1$. Therefore, $\varphi_{2,3}^{(u)}(t, y, 0) \equiv \varphi_{2,3}(t, y, 0)$ for every $u \in R_1$.

Analogously, the max process $\{\xi_0(t), t \geq 0\}$ and the sum process $\{\kappa_0(t), t \geq 0\}$ are independent if and only if the function $\Pi_3^{(u)}((w, \infty)) \equiv 0$ for $u > u_\pi$, $w > 0$. Now, $d^{(u)} \equiv d$ and the measure $\hat{\Pi}_3^{(u)}(C) \equiv \Pi_3(A)$ for every $u \in R_1$. Again, $\varphi_{2,3}^{(u)}(t, 0, z) \equiv \varphi_{2,3}(t, 0, z)$ for every $u \in R_1$. Obviously, the processes $\{\xi_0(t), t \geq 0\}$ and $\{\tau_0(t), t \geq 0\}$ are also independent.

In the case where both marginal independence assumptions are valid, the max process $\{\xi_0(t), t \geq 0\}$ and the two-dimensional sum process $\{(\gamma_0(t), \kappa_0(t)), t \geq 0\}$ are also independent. Indeed, $\Pi_{2,3}^{(u)}(((-\infty, -v] \cup (v, \infty)) \times (w, \infty)) \leq \Pi_2^{(u)}((-\infty, -v] \cup (v, \infty)) \wedge \Pi_3^{(u)}((w, \infty)) \equiv 0$ for every $u > u_\pi$, $v, w > 0$. Again, the constants satisfy $a^{(u)} \equiv a$, $d^{(u)} \equiv d$ while the measure $\hat{\Pi}_{2,3}^{(u)}(A) \equiv \Pi_{2,3}(A)$ for all $u \in R_1$. As a result, $\varphi_{2,3}^{(u)}(t, y, z) \equiv \varphi_{2,3}(t, y, z)$ for every $u \in R_1$.

4.2. *Identical components.* Let us consider the special case when the max and sum variables coincide, that is, $\xi_{\varepsilon,k} \equiv \gamma_{\varepsilon,k}, k = 1, 2, \ldots$. In this case, Condition A is implied by Condition B. Moreover, in this case, Condition B(a) implies that $\pi_1(u) = \Pi_2((u, \infty))$ for $u > 0$. Also $u_\pi = 0$ as follows from Condition B(a), which implies that $n_\varepsilon P\{\gamma_{\varepsilon,1} > u\} \to \infty$ as $\varepsilon \to 0$ for any $u < 0$. Furthermore, Condition C is implied by Condition B, and $\Pi_{2,3}^{(u)}((v, \infty) \times (w, \infty)) = \Pi_{2,3}((u \vee v, \infty) \times (w, \infty))$ for $u, v, w > 0$ while $\Pi_{2,3}^{(u)}((-\infty, -v] \times (w, \infty)) = 0$ for $u, v, w > 0$.

By summarizing all of the above information, we see that, in this case, the conditions to be validated are thoroughly reduced to Condition B.

Looking at the max processes, it is obvious that in this case $\xi_\varepsilon(t) = f_t(\gamma_\varepsilon(\cdot))$, where $f_t(x(\cdot)) = \max_{s \in (0,t]} \Delta_s(x(\cdot))$ and $\Delta_s(x(\cdot)) = x(s) - x(s-0)$ for $\varepsilon$ such that $n_\varepsilon t \geq 1$. So $\xi_\varepsilon(t)$ is a maximal jump of the process $\gamma_\varepsilon(s)$ on the interval $[0, t]$. This functional is a.s. $J$-continuous with respect to measure generated by the limiting process $\{\gamma_0(s), s \geq 0\}$ on the Borel $\sigma$-algebra of the space $D_1$ for every $t > 0$. Therefore, $\xi_0(t) = f_t(\gamma_0(\cdot)), t > 0$.

4.3. *Transformed max–sum processes with renewal stopping.* Let $f(t, x)$ be a continuous function which is defined on $[0, \infty) \times R_3$ while taking values in $R_1$. The transformed stochastic process $\{f(t, \zeta_\varepsilon(\tau_\varepsilon(t))), t > 0\}$ has trajectories which belong to the space $D_1$ with probability 1. Theorem 3.2 implies $J$-convergence of transformed processes

$$(4.1) \quad \{f(t, \zeta_\varepsilon(\tau_\varepsilon(t))), t \geq 0\} \overset{J}{\to} \{f(t, \zeta_0(\tau_0(t))), t > 0\} \quad \text{as } \varepsilon \to 0.$$

Below, several examples illustrate (4.1). We can apply the relationship to a variety of processes representing modifications of the original sum processes and max processes.

As a first example, take

$$\begin{aligned}
(4.2) \qquad & \{\gamma_\varepsilon(\tau_\varepsilon(t)) - \xi_\varepsilon(\tau_\varepsilon(t)), t > 0\} \\
& \xrightarrow{J} \{\gamma_0(\tau_0(t)) - \xi_0(\tau_0(t)), t > 0\} \qquad \text{as } \varepsilon \to 0.
\end{aligned}$$

As other examples, look at

$$(4.3) \qquad \left\{ \frac{\xi_\varepsilon(\tau_\varepsilon(t))}{a + |\gamma_\varepsilon(\tau_\varepsilon(t))|}, t > 0 \right\} \xrightarrow{J} \left\{ \frac{\xi_0(\tau_0(t))}{a + |\gamma_0(\tau_0(t))|}, t > 0 \right\} \qquad \text{as } \varepsilon \to 0$$

and

$$(4.4) \qquad \left\{ \frac{\xi_\varepsilon(\tau_\varepsilon(t))}{at + |\gamma_\varepsilon(\tau_\varepsilon(t))|}, t > 0 \right\} \xrightarrow{J} \left\{ \frac{\xi_0(\tau_0(t))}{at + |\gamma_0(\tau_0(t))|}, t > 0 \right\} \qquad \text{as } \varepsilon \to 0.$$

Here $a > 0$ is a regularization parameter which prevents the denominator in the last two relationships to take the value zero.

Relationships (4.2)–(4.4) establish the weak convergence of the functionals that describe deviations of max and sum processes with renewal stopping.

It is obvious that the process $\{\zeta_0(\tau_0(t)), t > 0\}$ and the transformed process $\{f(t, \zeta_0(\tau_0(t))), t > 0\}$ are stochastically continuous in points $t \in V$, where $V$ is the set of stochastic continuity of the internal stopping process $\{\tau_0(t), t > 0\}$. That is why (4.1) implies, for example, that for any $0 < T_1 < T_2 < \infty, T_1, T_2 \in V$,

$$(4.5) \qquad \sup_{t \in [T_1, T_2]} f(t, \zeta_\varepsilon(\tau_\varepsilon(t))) \implies \sup_{t \in [T_1, T_2]} f(t, \zeta_0(\tau_0(t))) \qquad \text{as } \varepsilon \to 0.$$

This relationship, applied to the modified processes (4.2)–(4.4), establishes the weak convergence of the functional that describes the maximal deviations given by the corresponding processes.

We stress that the results presented in the current article give the general framework for the asymptotic study of hybrid max–sum processes with renewal stopping in the form of weak and functional limit theorems. Questions about the explicit form of the corresponding limiting distributions require explicit calculations outside the scope of this article.

4.4. *Applications to risk processes.* Let $\{(\beta_{\varepsilon,n}, \kappa_{\varepsilon,n}), n = 1, 2, \ldots\}$ be, for every $\varepsilon > 0$, a sequence of i.i.d. random variables taking values in $R_1 \times [0, \infty)$. Further we need nonrandom functions $n_\varepsilon > 0$ for which $n_\varepsilon \to \infty$ as $\varepsilon \to 0$ and constants $c_\varepsilon \geq 0$. Let us introduce a process $\{\mu_\varepsilon(t) = c_\varepsilon t - \beta_\varepsilon(\tau_\varepsilon(t)), t \geq 0\}$, where $\beta_\varepsilon(t) = \sum_{k \leq tn_\varepsilon} \beta_{\varepsilon,k}, \tau_\varepsilon(t) = \sup(s : \kappa_\varepsilon(s) \leq t)$ and $\kappa_\varepsilon(t) = \sum_{k \leq tn_\varepsilon} \kappa_{\varepsilon,k}, t \geq 0$. Random variables $\beta_{\varepsilon,n}$ can be interpreted as successive claims within a portfolio, while $\kappa_{\varepsilon,n}$ stands for the interarrival times of the claims. The process $\tau_{\varepsilon}(t)$ counts the number of claims that arrived within the portfolio up to time $t$. The constant $c_{\varepsilon}$ acts as a premium rate. The resulting process $\{\mu_{\varepsilon}(t), t \geq 0\}$ is called the (risk) reserves process.

(a) This process can be represented in the form $\mu_{\varepsilon}(t) = c_{\varepsilon} \cdot (t - \kappa_{\varepsilon}(\tau_{\varepsilon}(t))) + \gamma_{\varepsilon}(\tau_{\varepsilon}(t)), t \geq 0$, where $\gamma_{\varepsilon}(t) = \sum_{k \leq t n_{\varepsilon}} (c_{\varepsilon} \kappa_{\varepsilon,k} - \beta_{\varepsilon,k}), t \geq 0$.

(b) This representation can be supplemented by the estimate $|c_{\varepsilon} \cdot (t - \kappa_{\varepsilon}(\tau_{\varepsilon}(t)))| \leq \xi_{\varepsilon}(\tau_{\varepsilon}(t)), t \geq 0$, where $\xi_{\varepsilon}(t) = \max_{k \leq 1 \vee t n_{\varepsilon}} c_{\varepsilon} \kappa_{\varepsilon,k}, t \geq 0$.

The representation (a) and the estimate (b) show that asymptotic behavior of the risk processes can be studied via mixed max–sum processes with renewal stopping based on the sequences of i.i.d. random variables $\{(\xi_{\varepsilon,n}, \gamma_{\varepsilon,n}, \kappa_{\varepsilon,n}), \; n = 1, 2, \ldots\}$, where $\xi_{\varepsilon,n} = c_{\varepsilon} \kappa_{\varepsilon,n}, \gamma_{\varepsilon,n} = c_{\varepsilon} \kappa_{\varepsilon,n} - \beta_{\varepsilon,n}, n = 1, 2, \ldots$.

(c) Assume that the function $\pi_1(u) = 0$ for $u > 0$ in the corresponding Condition A. Then automatically $u_\pi = 0$, since the random variables $\xi_{\varepsilon,n} = c_{\varepsilon} \kappa_{\varepsilon,n}$ are nonnegative.

Under Condition A and assumption (c), the corresponding limiting max process $\xi_0(t) \equiv 0, t > 0$. Thus, due to (b), the processes $c_{\varepsilon} \cdot (t - \kappa_{\varepsilon}(\tau_{\varepsilon}(t))) \xrightarrow{P} 0$ as $\varepsilon \to 0$ for $t > 0$. If Conditions B and C hold together with $D_1$ in the form $c > 0$ [this does not contradict (c)], then Theorem 3.2 and (b) imply that the processes $\{c_{\varepsilon} \cdot (t - \kappa_{\varepsilon}(\tau_{\varepsilon}(t))), t > 0\}$ $J$-converge as $\varepsilon \to 0$ to the process which is identically equal to zero. In this case, the asymptotic behavior of the risk processes $\{\mu_{\varepsilon}(t), t \geq 0\}$ coincides with the asymptotic behavior of the sum processes with renewal stopping $\{\gamma_{\varepsilon}(\tau_{\varepsilon}(t)), t \geq 0\}$. The latter is described in Theorems 3.1 and 3.2.

Many questions from current day actuarial science can be formulated in terms of functionals on the $\{(\beta_{\varepsilon,n}, \kappa_{\varepsilon,n}), \; n = 1, 2, \ldots\}$ process. The important concept of *large claim* can be formulated in terms of the portion that the maximal claim consumes of the total claim amount. The behavior of the ratio between the largest claim and the totality of the claims for the i.i.d. case can be derived from results in Breiman (1965). In Ammeter (1964) we find a first attempt to consider the effect of reducing the total claim amount by subtracting the maximal claim from it. An extension to measure the influence of the largest claims on their total is linked to the concept of Lorenz curve, which was already hinted at in Aebi, Embrechts and Mikosch (1992).

4.5. *More general stopping processes.* The method of proof used in Theorems 3.1 and 3.2 can be applied to more general positive nondecreasing stopping processes $\{\tau_{\varepsilon}(t), t > 0\}$ such that:

(i) For every $t > 0$, $\tau_\varepsilon(t)$ is a Markov moment for the max–sum processes $\{\hat{\zeta}_\varepsilon^{(h_0)}(t), t > 0\}$ for some $h_0 > -\infty$ if $u_\pi = -\infty$ and $h_0 = u_\pi$ if $u_\pi > -\infty$.

(j) The relationship of joint weak convergence (3.1) holds for processes $\{\zeta_\varepsilon(t), t > 0\}$ and $\{\tau_\varepsilon(t), t > 0\}$, with some set $V$ dense in $(0, \infty)$.

Note that we involve truncated processes $\{\hat{\zeta}_\varepsilon^{(h)}(t), t > 0\}$ only for one value $h = h_0$. Indeed, any Markov moment for the process $\{\zeta_\varepsilon^{(h_0)}(t), t > 0\}$ is also a Markov moment for the process $\{\hat{\zeta}_\varepsilon^{(h)}(t), t > 0\}$ for every $h \in [-\infty, h_0]$. Also, since $\max(x, h)$ is a continuous function in $x$, (j) implies that the same condition holds for processes $\{\hat{\zeta}_\varepsilon^{(h)}(t), t > 0\}$ and $\tau_\varepsilon(t), t > 0$, for any $h \in [-\infty, \infty)$.

Note also that in many cases one can avoid truncation in Condition (i). For example, if $\tau_\varepsilon(t)$ is, for every $t > 0$, a Markov moment for the max–sum processes $\{(\gamma_\varepsilon(t), \kappa_\varepsilon(t)), t > 0\}$, then condition (i) holds for any $h_0 \in R_1$.

A large number of various examples can be constructed where the stopping moments have a renewal structure. We refer to Silvestrov (1974, 2000), where one can find a variety of examples of stopping moments of such a type. Here we point out only two typical examples.

The first one is concerned with the case where $\{\tau_\varepsilon(t) = \sup(s : f(s, \hat{\zeta}_\varepsilon^{(h_0)}(s)) \leq t), t > 0\}$. Here $f(t, x)$ is a continuous function which is defined on $[0, \infty) \times R_3$ while taking values in $R_1$. By varying the function $f$, one can construct many types of renewal type stopping moments which are of a first-passage-time type.

To avoid consideration of improper moments, we assume that $\sup_{s \leq t} f(s, \hat{\zeta}_\varepsilon^{(h_0)}(s))$ tends to $\infty$ in probability as $t \to \infty$ for every $\varepsilon \in [0, \varepsilon_0)$. The alternative could be to truncate moments $\tau_\varepsilon(t)$. Obviously, condition (i) holds. It can be shown [see, e.g., Silvestrov (1974) or Whitt (1980)] that, under Conditions A–C, the functional $\sup(s : f(s, x(s)) \leq t)$ is for every $t > 0$ a.s. $J$-continuous with respect to the limiting process $\{\hat{\zeta}_0^{(h_0)}(s), s > 0\}$ if the process $\{\sup_{s \leq t} f(s, \hat{\zeta}_0^{(h_0)}(s)), t > 0\}$ is a.s. strictly monotonic. That is why condition (j) also holds with $V = (0, \infty)$.

Another example which we would like to point out covers the models with renewal extremal stopping. Such models were studied in Shanthikumar and Sumita (1983), Sumita and Shanthikumar (1985), Silvestrov and Teugels (1998b), Gut and Hüsler (1999) and Gut (2001). In this case $\{\tau_\varepsilon(t) = \sup(s : \Delta_s(\xi_\varepsilon(\cdot)) \leq t), t > 0\}$.

For simplicity, we assume that $\pi_1(u) > 0$ for all $u > 0$ to exclude the case of improper renewal moments for $\varepsilon$ small enough. Again condition (i) holds. It can be easily shown that, under Condition A, the functionals $\Delta_t(x(\cdot))$ and $\sup(s : \Delta_s(x(\cdot)) \leq t)$ are a.s. $J$-continuous with respect the limiting process $\{\xi_0(s), s \geq 0\}$ for every $t \in U$, where $U$ is the set of $u > 0$ which are points of continuity of the function $\pi_1(u)$ appearing in Condition A. So, condition (j) also holds.

Department of Mathematics and Physics
Mälardalen University
SE-721 23 Västerås
Sweden
e-mail: dmitrii.silvestrov@mdh.se

Department of Mathematics
Katholieke Universiteit Leuven
B-3001 Leuven (Heverlee)
Belgium
e-mail: jef.teugels@wis.kuleuven.ac.be